\documentclass[11pt]{article}
\usepackage{amsmath}
\usepackage{amsfonts}
\usepackage[bookmarksnumbered=true]{hyperref}

\newtheorem{theorem}{Theorem}[section]
\newtheorem{lemma}[theorem]{Lemma}
\newtheorem{definition}[theorem]{Definition}
\newtheorem{remark}[theorem]{Remark}

\begin{document}

\title{Backward Doubly Stochastic Differential Equations
with Jumps and Stochastic Partial Differential-Integral Equations
\thanks{This work is supported by National Natural Science Foundation of China Grant
10771122, Natural Science Foundation of Shandong Province of China
Grant Y2006A08 and National Basic Research Program of China (973
Program, No.2007CB814900)}}
\author{Qingfeng Zhu$^{\rm a}$ and Yufeng Shi$^{\rm b}$\thanks{Corresponding author, E-mail: yfshi@sdu.edu.cn}\\
{\small $^{\rm a}$ School of Statistics and Mathematics, Shandong
University of Finance},\\
{\small Jinan 250014, China}\\
{\small$^{\rm b}$School of Mathematics, Shandong University, Jinan 250100, China}}
\maketitle

\begin{abstract}In this paper, we study backward doubly stochastic differential equations driven
by Brownian motions and Poisson process (BDSDEP in short) with
non-Lipschitz coefficients on random time interval.
The probabilistic interpretation for the solutions to a class of
quasilinear stochastic partial differential-integral equations
(SPDIEs in short) is tthe solutionBDSDEP. Under non-Lipschitz
conditions, the existence and uniqueness results for measurable
solutions of BDSDEP are established via the smoothing technique.
Then, the continuous dependence for solutions of BDSDEP is derived. Finally,
the probabilistic interpretation for the solutions to a class of
quasilinear SPDIEs is given.\\
\indent{\it keywords:} Backward doubly stochastic differential equations,
 stochastic partial differential-integral equations, random measure, Poisson process
\end{abstract}

\section{Introduction}\label{sec:1}

Nonlinear backward stochastic differential equations with Brownian
motions as noise sources (BSDEs in short) have been independently
introduced by Pardoux and Peng \cite{PP1} and Duffie and Epstein \cite{DE}. By
virtue of BSDEs, Peng \cite{P} has given a probabilistic interpretation
(nonlinear Feynman-Kac formula) for the solutions of semilinear
parabolic partial differential equations (PDEs in short). In \cite{P},
Peng also gave an existence and uniqueness result of BSDEs with
random terminal time. And then Darling and Pardoux \cite{DP} proved an
existence and uniqueness result for BSDEs with random terminal time
under different assumptions. They applied their result to construct
a continuous viscosity solution for a class of semilinear elliptic
PDEs.

A class of backward doubly stochastic differential equations (BDSDEs
in short) was introduced by Pardoux and Peng \cite{PP2} in 1994, in order
to provide a probabilistic interpretation for the solutions of a
class of semilinear stochastic partial differential equations (SPDEs
in short). They have proved the existence and uniqueness of
solutions for BDSDEs under uniformly Lipschitz conditions. Since
then, Shi  et al. \cite{SGL} have relaxed the Lipschitz assumptions to
linear growth conditions. Bally and Matoussi \cite{BM} have given a
probabilistic interpretation of the solutions in Sobolev spaces for
semilinear parabolic SPDEs in terms of BDSDEs. Zhang and Zhao \cite{ZZ}
have proved the existence and uniqueness of solution for BDSDEs on
infinite horizons, and described the stationary solutions of SPDEs
by virtue of the solutions of BDSDEs on infinite horizons.

BSDEs driven by Brownian motions and Poisson process (BSDEP in
short) were first discussed by Tang and Li \cite{TL}. After then Situ
\cite{S} obtained an existence and uniqueness result for BSDEP with
non-Lipschitz coefficients, so as to give a probabilistic
interpretation for solutions of partial differential-integral
equations (PDIEs in short). Barles et al. \cite{BBP} and Yin and Mao \cite{YM}
discussed viscosity solutions for a system of partial
differential-integral equations in terms of BSDEs with jumps.
Recently BDSDEs driven by Brownian motions and Poisson process
(BDSDEP in short) with Lipschitzian coefficients on a fixed time
interval were discussed by Sun and Lu \cite{SL}.

Because of their important significance to SPDEs, it is necessary to
give intensive investigation to the theory of BDSDEs. In this paper
we study BDSDEs driven by Brownian motions and Poisson process
(BDSDEP in short) with non-Lipschitzian coefficients on random time
interval. Here the coefficients are assumed to be weaker than linear
growth, continuous and to satisfy some weak ``monotone" condition.
BDSDEP can provide more extensive frameworks for the probabilistic
interpretations (so-called nonlinear stochastic Feynman-Kac formula)
for the solutions of a class of quasilinear stochastic partial
differential-integral equations (SPDIEs in short). First, We
establish the existence and uniqueness results for measurable
solutions of BDSDEP based on the smoothing technique. Then we
discuss the continuous dependence for solutions of BDSDEP.
Finally, by virtue of BDSDEP, we show the probabilistic
interpretation for the solutions of a class of quasilinear SPDIEs.

The paper is organized as follows. In Section  \ref{sec:2}, the basic assumptions
are given. In Section  \ref{sec:3}, the existence and
uniqueness for BDSDEP with non-Lipschitz coefficients on
random time interval is proved. In Section  \ref{sec:4}, the continuous
dependence for solutions of BDSDEP  is discussed.
Finally, in Section  \ref{sec:5}, the probabilistic interpretation
for the solutions to a class of quasilinear SPDIEs is
given by virtue of this class of BDSDEP.

\section{Setting of the problem}\label{sec:2}

Let $(\Omega,{\cal F},P) $ be a complete probability space, and
$[0,T]$ be an arbitrarily large fixed time duration throughout this
paper. We suppose $\{{\cal F}_t\}_{t \geq 0}$ is generated by the
following three mutually independent processes:

(i)\hspace{0.1cm}Let $\left\{ W_t;0\leq t\leq T\right\}$ and
$\left\{ B_t;0\leq t\leq T\right\} $ be two standard Brownian
motions defined on $( \Omega ,{\cal F},P) $, with values
respectively in $\mathbb{R}^d$ and in $\mathbb{R}^l$.

(ii)\hspace{0.1cm}Let $N$ be a Poisson random measure, on
$\mathbb{R}_{+} \times Z$, where $Z \subset \mathbb{R}^r$ is a
nonempty open set equipped with its Borel field ${\cal B}(Z)$, with
compensator $\widehat{N}({\rm d}z{\rm d}t)=\lambda({\rm d}z){\rm
d}t$, such that $\widetilde{N}(A \times [0,t])=(N-\widehat{N})(A
\times [0,t])_{t\geq 0}$ is a martingale for all $A \in {\cal B}(Z)$
satisfying $\lambda(A) < \infty$. $\lambda$ is assumed to be a
$\sigma$-finite measure on $(Z,{\cal B}(Z))$ and is called the
characteristic measure.

Let ${\cal N}$ denote the class of $P$-null elements of ${\cal F}$.
For each $t\in \left[ 0,T\right] $, we define $ {\cal F}_t\doteq
{\cal F}_t^W\vee {\cal F}_{t,T}^B\vee {\cal F}_t^N,$ where for any
process $\{\eta_t\},{\cal F}_{s,t}^\eta=\sigma
\left\{\eta_r-\eta_s;s\leq r\leq t\right\}\vee {\cal N}$, ${\cal
F}_t^\eta={\cal F}_{0,t}^\eta $. Note that the collection $\{
{\cal F}_t,t\in$ $\left[ 0,T\right] \} $ is neither increasing
nor decreasing, and it does not constitute a classical filtration.

Let $\tau=\{\tau(\omega)\}$ be an ${\cal F}_t$-measurable time on
$\left[ 0,T\right] $. We introduce the following notations:
\begin{eqnarray*}
S^2\left([0,\tau];\mathbb{R}^n\right)&=&\{v_t,0\leq t\leq \tau,\
 \mbox {is an}\ \mathbb{R}^n\mbox{-valued},\ {\cal F}_t \mbox{-measurable process }\\
  &&\mbox{such that} \ E(\sup_{0\leq t\leq\tau}|v_{t}|^{2})<\infty\},\\
M^2(0,\tau;\mathbb{R}^n)&=&\{v_t,0\leq t\leq \tau,\
 \mbox {is an}\ \mathbb{R}^n\mbox{-valued},\ {\cal F}_t\mbox{-measurable process }\\
 && \mbox{such that} \ E\int_{0}^\tau|v_{t}|^{2}dt<\infty\},\\
F^2_N(0,\tau;\mathbb{R}^n)&=&\{k_t,0\leq t\leq \tau,\
 \mbox {is an}\ \mathbb{R}^n\mbox{-valued},\ {\cal F}_t\mbox{-measurable process }\\
 && \mbox{such that} \ E\int_{0}^\tau\int_{Z}|k_{t}(z)|^{2}\lambda(dz)dt<\infty\},\\
L^2_{\lambda(\cdot)}(\mathbb{R}^n)&=&\{k(z),\ k(z)\ \mbox{is an}\
\mathbb{R}^n\mbox{-valued}, \ {\cal B}(Z)\mbox{-measurable function}\\
&&\mbox{ such that}\ \|k\|=(\int_Z|k(z)|^2\lambda(dz))^{1/2}<\infty\},\\
L^{2}(\Omega, {\cal F}_\tau,P;\mathbb{R}^{n})&=&\{\xi,\ \xi\
\mbox{is an}\ \mathbb{R}^n\mbox{-valued},\  {\cal
F}_\tau\mbox{-measurable  random variable}\\
 && \mbox{such that}\ E|\xi|^2<\infty\}.
\end{eqnarray*}

Consider the following BDSDE with Brownian motions and Poisson
Process (BDSDEP in short):

\begin{eqnarray}\label{eq:1}
\nonumber P_t&=&\xi+\int_{t\wedge\tau}^\tau f(s,P_s,Q_s,K_s){\rm
d}s+\int_{t\wedge\tau}^\tau g(s,P_s,Q_s,K_s)dB_s\\
&&-\int_{t\wedge\tau}^\tau Q_sdW_s-\int_{t\wedge\tau}^\tau
\int_ZK_s(z)\widetilde{N}(dzds),\quad t\geq 0,
\end{eqnarray}
where $\xi\in L^{2}(\Omega, {\cal F}_\tau,P;\mathbb{R}^{n})$,
\begin{eqnarray*}
f: \Omega \times [0,T] \times \mathbb{R}^{n} \times \mathbb{R}^{n
\times d} \times L^2_{\lambda(\cdot)}(\mathbb{R}^n)\rightarrow
\mathbb{R}^{n},\end{eqnarray*}
and
\begin{eqnarray*}g: \Omega \times [0,T] \times
\mathbb{R}^{n} \times \mathbb{R}^{n \times d} \times
L^2_{\lambda(\cdot)}(\mathbb{R}^n) \rightarrow \mathbb{R}^{n \times
l}.
\end{eqnarray*}
We note that the integral with respect to $\{B_t\}$ is a ``backward
It\^o integral" and the integral with respect to $\{W_t\}$ is a
standard forward It\^o integral. These two types of integrals are
particular cases of the It\^o-Skorohod integral (see Pardoux and
Peng \cite{PP2}). We use the usual inner product $\langle \cdot ,\cdot
\rangle $ and Euclidean norm $ | \cdot | $ in $\mathbb{R}^n$,
$\mathbb{R}^{n\times l}$ and $\mathbb{R}^{n\times d}.$ All the
equalities and inequalities mentioned in this paper are in the sense
of $dt\times dP$ almost surely on $\left[ 0,\tau\right] \times
\Omega $.

\begin{definition}\label{def:2.1}
A solution of BDSDEP (\ref{eq:1}) is a triple of ${\cal F}_t$-measurable
stochastic processes $(P,Q,K)$ which
belongs to the space $ S^2([0,\tau];\mathbb{R}^{n})\times
M^2(0,\tau$; $\mathbb{R}^{n\times d})\times
F^{2}_{N}(0,\tau;\mathbb{R}^{n})$ and satisfies BDSDEP (\ref{eq:1}).
\end{definition}

We assume that
\begin{enumerate}
\item[(H1)]
$\xi\in L^{2}(\Omega, {\cal F}_\tau,P;\mathbb{R}^{n})$;
\item[(H2)]
$f(t,p,q,k)$, $g(t,p,q,k)$ are continuous in $(p,q,k)\in
\mathbb{R}^{n} \times \mathbb{R}^{n \times d} \times L^2_{\lambda(\cdot)}(\mathbb{R}^n)$;
\item[(H3)]
 $f=f_1+f_2$, $f_i=f_i(t,p,q,k):\Omega \times [0,T] \times \mathbb{R}^{n} \times \mathbb{R}^{n
\times d} \times L^2_{\lambda(\cdot)}(\mathbb{R}^n)\rightarrow
\mathbb{R}^{n}$, $i=1,2$, and $g(t,p,q,k)$ are ${\cal
F}_t$-measurable processes, such that for all $t\in[0,T]$;
$p,p_1,p_2\in\mathbb{R}^{n}$; $q,q_1,q_2\in\mathbb{R}^{n\times d}$;
$k,k_1,k_2\in L^2_{\lambda(\cdot)}(\mathbb{R}^n)$,
\begin{eqnarray*}
&&|f_1(t,p,q,k)|\leq \mu(t),\\
&&|f_2(t,p,q,k)|\leq \mu(t)(1+|p|+|q|+\|k\|),\\
&&|g(t,p,q,k)|\leq \mu(t),
\end{eqnarray*}
where $\mu(t)\geq 0$ is real and non-random function such that
$\bar \mu=\int_0^T\mu^2(t)dt< \infty.$
\item[(H4)]
for all $t\in[0,T]$; $p,p_1,p_2\in\mathbb{R}^{n}$;
$q,q_1,q_2\in\mathbb{R}^{n\times d}$; $k,k_1,k_2\in
L^2_{\lambda(\cdot)}(\mathbb{R}^n)$, such that
\begin{eqnarray*}
&
&\left<p_1-p_2,f_1(t,p_1,q_1,k_1)-f_1(t,p_2,q_2,k_2)\right>\\
&&\leq \mu(\rho(|p_1-p_2|^2)+|p_1-p_2|(|q_1-q_2|+\|k_1-k_2|)),\\
&& |f_1(t,p,q,k_1)-f_1(t,p,q,k_2)|\leq \mu\|k_1-k_2\|,\\
&&|f_2(t,p_1,q_1,k_1)-f_2(t,p_2,q_2,k_2)|\\
&&\leq \mu(|p_1-p_2|+|q_1-q_2|+\|k_1-k_2\|),\\
&&|g(t,p_1,q_1,k_1)-g(t,p_2,q_2,k_2)|^2 \\
&&\leq \mu(|p_1-p_2|^2+|p_1-p_2|(|q_1-q_2|+\|k_1-k_2|)),
\end{eqnarray*}
where $\mu >0$ is a constant, and $\rho(\cdot)$
is a  nondecreasing, continuous and concave function from $R_+$
to $R_+$ such that $\rho(0)=0$, $\rho(u)>0$, as $u>0$, and
$\int_{0^+}du/\rho(u)=+\infty. $
\end{enumerate}

\section{Existence and uniqueness of solutions to BDSDEP with non-Lipschitz coefficients}
\label{sec:3}

In order to prove the existence and uniqueness results of solutions to BDSDEP
with non-Lipschitz coefficients on random time interval, we introduce the following
lemmas and theorems.

\begin{lemma}\label{lem:3.1} (A priori estimate).
Under the assumption (H3). If $(P_t, Q_t, K_t)$ is a solution of (\ref{eq:1}), then
\begin{eqnarray*}
E\left(\sup\limits_{t\leq \tau}|P_t|^2+\int_0^\tau|Q_t|^2dt
+\int_{0}^{\tau}\|K_{t}\|^{2}dt\right)\leq C_T<\infty,
\end{eqnarray*}
where $C_T \geq 0$ is a constant depending on $T$, $\int_0^T
\mu^2(t)dt$ and $E|\xi|^2$ only.
\end{lemma}

\noindent{\bf Proof.} From (H3), we easily have
\begin{eqnarray*}
\left<p,f(t,p,q,k)\right>\leq \mu(t)(1+2|p|^2+|p|(|q|+\|k\|)),
\end{eqnarray*}
where $\mu(t)$ has the property stated in (H3).
Applying It\^o's formula to $|P_t|^2$, we have
\begin{eqnarray*}
&&E\left(|P_{t\wedge\tau}|^2+\int_{t\wedge\tau}^\tau|Q_s|^2ds
+\int_{t\wedge\tau}^\tau\|K_{s}\|^{2}ds\right) \\
&=&E|\xi|^2+2E\int_{t\wedge\tau}^\tau\left<P_s,f(s,P_s,Q_s,K_s)\right>ds
+E\int_{t\wedge\tau}^\tau|g(s,P_s,Q_s,K_s)|^2ds\\
&\leq&
E|\xi|^2+ 2E\int_{t\wedge\tau}^\tau
\mu(s)(1+2|P_s|^2+|P_s|(|Q_s|+\|K_{s}\|))
ds+\bar \mu,
\end{eqnarray*}
we deduce
\begin{eqnarray*}
&&E\left(|P_{t\wedge\tau}|^2+\frac{1}{2}\int_{t\wedge\tau}^\tau|Q_s|^2ds
+\frac{1}{2}\int_{t\wedge\tau}^\tau\|K_{s}\|^{2}ds\right)\\
&\leq&E|\xi|^2+T+2\bar \mu+
E\int_{t}^T(4\mu(s)+2\mu^2(s))|P_s|^2ds.
\end{eqnarray*}
By Gronwall inequality, we have
\begin{eqnarray*}
E\left(|P_{t\wedge\tau}|^2+\frac{1}{2}\int_{t\wedge\tau}^\tau|Q_s|^2ds
+\frac{1}{2}\int_{t\wedge\tau}^\tau\|K_{s}\|^{2}ds\right) \leq
\widetilde C_T,
\end{eqnarray*}
where
\begin{eqnarray*}
 \widetilde C_T=\left(E|\xi|^2+T+2\bar \mu\right)
 \exp{\left(\int_{0}^T(4\mu(s)+2\mu^2(s))ds\right)}.
\end{eqnarray*}
In particular,
\begin{eqnarray*}
E\left(|P_{0}|^2+\frac{1}{2}\int_{0}^\tau|Q_s|^2ds
+\frac{1}{2}\int_{0}^\tau\|K_{s}\|^{2}ds\right) \leq \widetilde C_T.
\end{eqnarray*}

Applying It\^o's formula to $|P_t|^2$ on $[0,t\wedge \tau]$, we have
\begin{eqnarray*}
&&|P_{t\wedge \tau}|^2\\
&=& |P_0|^2+2\int_0^{t\wedge \tau}\langle
P_s,f(s,P_s,Q_s,K_s)\rangle ds+2\int_0^{t\wedge
\tau}\langle P_s,g(s,P_s,Q_s,K_s)\rangle dB_s\\
&&-2\int_0^{t\wedge\tau}\langle P_s,Q_s\rangle
dW_s+2\int_0^{t\wedge\tau}\int_Z\langle
P_s,K_s(z)\rangle\widetilde{N}(dzds)\\
&&+\int_0^{t\wedge\tau}|g(s,P_s,Q_s,K_s)|^2ds
-\int_0^{t\wedge\tau}|Q_s|^2ds +\int_0^{t\wedge\tau}\|K_s\|^2ds.
\end{eqnarray*}
Taking supremum and expectation, we get
\begin{eqnarray*}
&&E\sup\limits_{t \leq \tau}|P_{t\wedge\tau}|^2\\
&\leq& E|P_0|^2+2E\int_0^{\tau}\mu(s)(1+2|P_s|^2+|P_s|(|Q_s|+\|K_s\|))
ds+\int_0^{T}\mu^2(s)ds\\
&& +E\int_0^{\tau}|Q_s|^2ds
+E\int_0^{\tau}\|K_s\|^2ds+2E\sup\limits_{t\leq
\tau}|\int_0^{t\wedge\tau}\langle
P_s,g(s,P_s,Q_s,K_s)\rangle dB_s|\\
&&+2E\sup\limits_{t\leq \tau}|\int_0^{t\wedge\tau}\langle
P_s,Q_s\rangle dW_s|+2E\sup\limits_{t\leq
\tau}|\int_0^{t\wedge\tau}\int_Z\langle
P_s,K_s(z)\rangle\widetilde{N}(dzds)|.
\end{eqnarray*}
By Burkholder-Davis-Gundy's inequality, we deduce
\begin{eqnarray*}
&&E\left(\sup_{t\leq \tau}\left|\int_0^{t\wedge\tau}\langle
P_s,g(s,P_s,Q_s,K_s)\rangle dB_s\right|\right)\\
& \leq &  c E\left(\int_0^\tau|P_{s\wedge\tau}|^2\cdot |g(s,P_s,Q_s,K_s)|^2ds\right)^{1/2}\\
& \leq & c E\left(\left(\sup_{t\leq \tau
}|P_{t\wedge\tau}|^2\right)^{1/2}
\left(\int_0^\tau|g(s,P_s,Q_s,K_s)|^2ds\right)^{1/2}\right)\\
& \leq &\frac18E\sup_{t\leq \tau }|P_{t\wedge\tau}|^2
+2c^2\int_0^{T}\mu^2(s)ds.
\end{eqnarray*}
In the same way, we have
\begin{eqnarray*}
&&E\sup\limits_{t\leq \tau}|\int_0^{t\wedge\tau}\langle
P_s,Q_s\rangle dW_s| \leq\frac18E\sup_{t\leq \tau
}|P_{t\wedge\tau}|^2  +2c^2\int_0^{\tau}|Q_s|^2ds,\\
&&E\sup\limits_{t\leq \tau}|\int_0^{t\wedge\tau}\int_Z\langle
P_s,K_s(z)\rangle\widetilde{N}(dzds)|\leq\frac18E\sup_{t\leq \tau
}|P_{t\wedge\tau}|^2 +2c^2\int_0^{\tau}\|K_s\|^2ds.
\end{eqnarray*}
Hence
\begin{eqnarray}
\nonumber E\sup\limits_{t \leq \tau}|P_{t\wedge\tau}|^2 &\leq&
4E|P_0|^2+4T+8\bar \mu+4E\int_0^{\tau}(4\mu(s)+2\mu^2(s))|P_{s\wedge\tau}|^2ds\\
\nonumber
&&+4(1+2c^2)E\int_0^{\tau}(|Q_s|^2+\|K_s\|^2)
ds \leq C_T<\infty.
\end{eqnarray}
\quad$\Box$

As a preparation for the study of BDSDEP (\ref{eq:1}),
we first discuss a simpler BDSDEP as follows
\begin{eqnarray}\label{eq:2}
\nonumber P_t&=&\xi+\int_{t\wedge\tau}^\tau f(s)ds+\int_{t\wedge\tau}^\tau
g(s)dB_s-\int_{t\wedge\tau}^\tau Q_s dW_s\\
&&-\int_{t\wedge\tau}^\tau
\int_ZK_s(z)\widetilde{N}(dzds),\quad t\geq 0.
\end{eqnarray}
We have
\begin{lemma}\label{lem:3.2}
Given $\xi\in L^{2}(\Omega, {\cal F}_\tau,P;\mathbb{R}^{n})$,
$f(t)\in M^2(0,\tau;\mathbb{R}^n)$ and
$g(t)\in M^2(0$, $\tau; \mathbb{R}^{n\times l})$,
then (\ref{eq:2}) has a unique solution in $ S^2([0,\tau];\mathbb{R}^{n})\times
M^2(0,\tau$; $\mathbb{R}^{n\times d})\times
F^{2}_{N}(0,\tau;\mathbb{R}^{n})$.
\end{lemma}
{\bf Proof.} {\it Uniqueness.}
Let $(P^1,Q^1)$ and $(P^2,Q^2)$ be two solution of (\ref{eq:2}).
Applying It\^o's formula to $|P^1_t-P^2_t|^2$, we have
\begin{eqnarray*}
E|P^1_{t\wedge\tau}-P^2_{t\wedge\tau}|^2+E\int_{t\wedge\tau}^\tau
|Q^1_s-Q^2_s|^2ds+E\int_{t\wedge\tau}^\tau \|K^1_s-K^2_s\|^2ds=0.
\end{eqnarray*}
Then
\begin{eqnarray*}
E|P^1_{t\wedge\tau}-P^2_{t\wedge\tau}|^2=0,\
E\int_{t\wedge\tau}^\tau|Q^1_s-Q^2_s|^2ds=0,
\end{eqnarray*}
and
\begin{eqnarray*}
E\int_{t\wedge\tau}^\tau\|K^1_s-K^2_s\|^2ds=0,\ 0\leq t \leq T.
\end{eqnarray*}
Hence $P^1_t=P^2_t$, $Q^1_t=Q^2_t$ and  $K^1_t=K^2_t$ a.s.. The
uniqueness is obtained.

{\it Existence.} We define the filtration $({\cal G}_t)_{0 \leq t
\leq T}$ by
\begin{eqnarray*}
{\cal G}_t={\cal F}^W_t \vee {\cal F}^B_T \vee {\cal F}^N_t
\end{eqnarray*}
and the ${\cal G}_t$-square integrable martingale
\begin{eqnarray*}
M_t=E\left[\xi + \int_0^\tau f(s)ds+\int_0^\tau g(s)dB_s|{\cal
G}_{t}\right],\quad t \geq 0.
\end{eqnarray*}
An obvious extension of It\^{o}'s martingale representation theorem
(see \cite{IW}) yields the existence of $(Q_s,K_s)$ such that
\begin{eqnarray*}
M_t = M_0 + \int_0^tQ_sdW_s + \int_0^t\int_Z
K_{s}(z)\widetilde{N}(dzds),
\end{eqnarray*}
and
\begin{eqnarray*}
E\int_0^T(|Q_s|^2+\|K_s\|^2)ds<\infty.
\end{eqnarray*}
Particularly
\begin{eqnarray*}
&&M_0 + \int_0^\tau Q_sdW_s + \int_0^\tau\int_Z
K_{s}(z)\widetilde{N}(dzds)\\
&=&M_{\tau}=\xi + \int_0^\tau f(s)ds+\int_0^\tau g(s)dB_s,
\end{eqnarray*}
or
\begin{eqnarray}\label{eq:3}
\nonumber&&M_0 + \int_0^{t\wedge\tau} Q_sdW_s + \int_0^{t\wedge\tau}\int_Z
K_{s}(z)\widetilde{N}(dzds)\\
\nonumber&=&\xi + \int_0^\tau f(s)ds+\int_0^\tau g(s)dB_s
-\int_{t\wedge\tau}^\tau Q_szdW_s- \int_{t\wedge\tau}^\tau\int_Z
K_{s}(z)\widetilde{N}(dzds).\\
\end{eqnarray}
We set
\begin{eqnarray*}
P_t=E\left[\xi + \int_{t\wedge\tau}^\tau f(s)ds+\int_{t\wedge\tau}^\tau g(s)dB_s|{\cal
G}_{t}\right],\quad t \geq 0.
\end{eqnarray*}
It follows that
\begin{eqnarray*}
P_t&=&E\left[\xi + \int_{0}^\tau f(s)ds+\int_{0}^\tau g(s)dB_s|{\cal
G}_{t}\right]\\
&&-E\left[\int_{0}^{t\wedge\tau} f(s)ds+\int_{0}^{t\wedge\tau} g(s)dB_s|{\cal
G}_{t}\right]\\
&=&M_0 + \int_0^{t\wedge\tau} Q_sdW_s + \int_0^{t\wedge\tau}\int_Z
K_{s}(z)\widetilde{N}(dzds)\\
&&-\int_0^{t\wedge\tau}f(s)ds-\int_0^{t\wedge\tau}g(s)dB_s.
\end{eqnarray*}
This with (\ref{eq:3}) implies that $(P_t,Q_t,K_t)$ solves (\ref{eq:2}). \quad $\Box$

In the following of this section we derive the existence and uniqueness results
for solutions of BDSDEP on random time interval with Lipschitzian and
non-Lipschitzian coefficients. The first one, that is Theorem \ref{thm:3.3},
deal with the case where $f$ is Lipschitz continuous.

\begin{theorem}\label{thm:3.3}
Under the assumptions (H1)-(H4), if $f_1=0$, (\ref{eq:1})has a unique solution
$(P_t, Q_t, K_t)$ in $ S^2([0,\tau];\mathbb{R}^{n})\times
M^2(0,\tau$; $\mathbb{R}^{n\times d})\times
F^{2}_{N}(0,\tau;\mathbb{R}^{n})$.
\end{theorem}

\noindent{\bf Proof.}\ We define recursively a sequence $\{(P_t^i,Q_t^i,K_t^i)\}_{i=0,1,\dots}$ as follows.
Let $P_t^0=0,Q_t^0=0,K_t^0=0$.  By Lemma \ref{lem:3.2}, for any $(P^i_t,Q^i_t,K^i_t)\in
S^2([0,\tau];\mathbb{R}^{n})\times
M^2(0,\tau$; $\mathbb{R}^{n\times d})\times
F^{2}_{N}(0,\tau;\mathbb{R}^{n})$, $i=1,2$ there exists unique $(P^{i+1}_t,Q^{i+1}_t,K^{i+1}_t)$,
satisfying
\begin{eqnarray*}
P^{i+1}_t&=&\xi+\int_{t\wedge\tau}^\tau f(s,\bar P^i_s,\bar Q^i_s,\bar
K^i_s) ds+\int_{t\wedge\tau}^\tau g(s,\bar P^i_s,\bar Q^i_s,\bar
K^i_s) dB_s\\
&&-\int_{t\wedge\tau}^\tau
Q^{i+1}_sdW_s-\int_{t\wedge\tau}^\tau
\int_ZK^{i+1}_s(z)\widetilde{N}(dzds),\ t\geq 0.
\end{eqnarray*}
Moreover, by Lemma \ref{lem:3.2}, $(P^{i+1}_t,Q^{i+1}_t,K^{i+1}_t)\in  S^2([0,\tau];\mathbb{R}^{n})\times
M^2(0,\tau$; $\mathbb{R}^{n\times d})\times
F^{2}_{N}(0,\tau;\mathbb{R}^{n})$.

Let $\bar P^{i+1}_t=P^{i+1}_t-P^{i}_t$, $\bar Q^{i+1}_t=Q^{i+1}_t-Q^{i}_t$, $\bar K^{i+1}_t=K^{i+1}_t-K^{i}_t$,
By the It\^{o}'s formula to $|P^{i+1}_{t\wedge\tau}-P^{i+1}_{t\wedge\tau}|^2e^{-\beta t}$, we have
\begin{eqnarray*}
&&E|P^1_{t\wedge\tau}-P^2_{t\wedge\tau}|^2e^{-\beta t}
+\beta E\int_{t\wedge\tau}^\tau|P^1_s-P^2_s|^2e^{-\beta s}ds\\
&&+E\int_{t\wedge\tau}^\tau|Q^1_s-Q^2_s|^2e^{-\beta s}ds
+E\int_{t\wedge\tau}^\tau\|K^1_s-K^2_s\|^2e^{-\beta s}ds\\
&=&2E\int_{t\wedge\tau}^\tau\left<P^1_s-P^2_s,
f(s,\bar P^1_s,\bar Q^1_s,\bar K^1_s)-f(s,\bar P^2_s,\bar Q^2_s,\bar K^2_s)\right>e^{-\beta s}ds\\
&&+E\int_{t\wedge\tau}^\tau|g(s,\bar P^1_s,\bar Q^1_s,\bar K^1_s)-g(s,\bar P^2_s,\bar Q^2_s,\bar K^2_s)|^2e^{-\beta s}ds
\end{eqnarray*}
\begin{eqnarray*}
&\leq&2\mu E\int_{t\wedge\tau}^\tau
|P_s^1-P_s^2|(|\bar P_s^1-\bar P_s^2|+|\bar Q_s^1-\bar Q_s^2|+\|\bar K_s^1-\bar K_s^2\|)e^{-\beta s}
ds\\
&&+\mu E\int_{t\wedge\tau}^\tau(|\bar P_s^1-\bar P_s^2|^2+|\bar P_s^1-\bar P_s^2|
(|\bar Q_s^1-\bar Q_s^2|+\|\bar K_s^1-\bar K_s^2\|))e^{-\beta s}ds\\
&\leq&\frac{1}{4}E\int_{t\wedge\tau}^\tau\left(|\bar P^1_s-\bar
P^2_s|^2+|\bar Q^1_s-\bar Q^2_s|^2+\|\bar K^1_s-\bar K^2_s\|^2 \right)e^{-\beta s}ds\\
&&+12\mu E\int_{t\wedge\tau}^\tau|P^1_s-P^2_s|^2e^{-\beta s}ds
+3\mu E\int_{t\wedge\tau}^\tau|\bar P^1_s-\bar P^2_s|^2e^{-\beta s}ds\\
&&+\frac{1}{4}E\int_{t\wedge\tau}^\tau\left(|\bar Q^1_s-\bar Q^2_s|^2
+\|\bar K^1_s-\bar K^2_s\|^2 \right)e^{-\beta s}ds,
\end{eqnarray*}
we deduce
\begin{eqnarray*}
&&E|P^1_{t\wedge\tau}-P^2_{t\wedge\tau}|^2e^{-\beta t}
+(\beta-12\mu)E\int_{t\wedge\tau}^\tau|P^1_s-P^2_s|^2e^{-\beta s}ds\\
&&+E\int_{t\wedge\tau}^\tau|Q^1_s-Q^2_s|^2e^{-\beta s}ds
+E\int_{t\wedge\tau}^\tau\|K^1_s-K^2_s\|^2e^{-\beta s}ds \\
&\leq&\frac{1}{2}E\int_{t\wedge\tau}^\tau\left(|\bar Q^1_s-\bar Q^2_s|^2
+\|\bar K^1_s-\bar K^2_s\|^2 \right)e^{-\beta s}ds\\
&&+(\frac{1}{4}+3\mu) E\int_{t\wedge\tau}^\tau|\bar P^1_s-\bar P^2_s|^2e^{-\beta s}ds,
\end{eqnarray*}
Now choose $\beta=12\mu+\displaystyle\frac{1+12\mu}{2}$,
and define $\bar c=\displaystyle\frac{1+12\mu}{2}$.
\begin{eqnarray*}
&&E|P^1_{t\wedge\tau}-P^2_{t\wedge\tau}|^2e^{-\beta t}\\
&&+\int_{t\wedge\tau}^\tau\left(\bar c|P^1_s-P^2_s|^2+|Q^1_s-Q^2_s|^2+\|K^1_s-K^2_s\|^2\right)e^{-\beta s}ds\\
&\leq&\frac{1}{2}E\int_{t\wedge\tau}^\tau\left(\bar c|\bar P^1_s-\bar P^2_s|^2+|\bar Q^1_s-\bar Q^2_s|^2
+\|\bar K^1_s-\bar K^2_s\|^2 \right)e^{-\beta s}ds,
\end{eqnarray*}
It follows immediately that
\begin{eqnarray*}
&&E\int_{t\wedge\tau}^\tau\left(\bar c|P^1_s-P^2_s|^2+|Q^1_s-Q^2_s|^2+\|K^1_s-K^2_s\|^2\right)e^{-\beta s}ds\\
&\leq&\frac{1}{2}E\int_{t\wedge\tau}^\tau\left(\bar c|\bar P^1_s-\bar P^2_s|^2+|\bar Q^1_s-\bar Q^2_s|^2
+\|\bar K^1_s-\bar K^2_s\|^2 \right)e^{-\beta s}ds,
\end{eqnarray*}
and $\{(P_t^i,Q_t^i,K_t^i)\}_{i=0,1,\dots}$ is a Cauchy sequence in
$S^2([0,\tau];\mathbb{R}^{n})\times
M^2(0,\tau$; $\mathbb{R}^{n\times d})\times
F^{2}_{N}(0,\tau;\mathbb{R}^{n})$, and that
\begin{eqnarray*}
\{(P_t,Q_t,K_t)\}=\lim\limits_{i\to \infty}\{(P_t^i,Q_t^i,K_t^i)\}
\end{eqnarray*}
solves (\ref{eq:1}). \quad $\Box$

The next theorem is the main result of this section, which generalizes
the result of Theorem \ref{thm:3.3} to the case  where $f$ is continuous but not
Lipschitz continuous.
\begin{theorem}\label{thm:3.4}
Under the assumptions (H1)-(H4), then (\ref{eq:1})has a unique solution
$(P_t,Q_t,K_t)$.
\end{theorem}
{\bf Proof.}{\it Uniqueness.}\
Let $(P^1_s,Q^1_s,K^1_s)$ and $(P^2_s,Q^2_s,K^2_s)$ be two solutions
of (\ref{eq:1}). Applying It\^{o}'s formula to $|P_s^1-P_s^2|^2$, we obtain
\begin{eqnarray*}
&&E\left(|P^1_{t\wedge\tau}-P^2_{t\wedge\tau}|^2+\int_{t\wedge\tau}^\tau|Q^1_s-Q^2_s|^2
ds +\int_{t\wedge\tau}^\tau\|K^1_s-K^2_s\|^2ds\right) \\
&=&2E\int_{t\wedge\tau}^\tau\left<P^1_s-P^2_s,
f(s,P^1_s,Q^1_s,K^1_s)-f(s,P^2_s,Q^2_s,K^2_s)\right>
ds\\
&&+E\int_{t\wedge\tau}^\tau|g(s,P^1_s,Q^1_s,K^1_s)-g(s,P^2_s,Q^2_s,K^2_s)|^2ds\\
&\leq&2\mu E\int_{t\wedge\tau}^\tau
(\rho(|P_s^1-P_s^2|^2)+|P_s^1-P_s^2|(|Q_s^1-Q_s^2|+\|K_s^1-K_s^2\|))
ds\\
&&+\mu E\int_{t\wedge\tau}^\tau(|P_s^1-P_s^2|^2+|P_s^1-P_s^2|(|Q_s^1-Q_s^2|+\|K_s^1-K_s^2\|))ds.
\end{eqnarray*}
From (H4), we have
\begin{eqnarray*}
X_t&=&E\left(|P^1_{t\wedge\tau}-P^2_{t\wedge\tau}|^2
+\frac{1}{2}\int_{t\wedge\tau}^\tau|Q^1_s-Q^2_s|^2ds
+\frac{1}{2}\int_{t\wedge\tau}^\tau\|K^1_s-K^2_s\|^2
ds\right)\\
&\leq&\mu E\int_{t\wedge\tau}^\tau(2\rho(|P_s^1-P_s^2|^2)
+11|P_s^1-P_s^2|^2)ds\\
&\leq& 11\mu\int_t^T\rho_1\left(X_s\right)ds,
\end{eqnarray*}
where
\begin{eqnarray*}
 \rho_1(u)=2\rho(u)+16u.
\end{eqnarray*}
By the Bahari's inequality, we obtain
\begin{eqnarray*}
E\left(|P^1_{t\wedge\tau}-P^2_{t\wedge\tau}|^2+\int_{t\wedge\tau}^\tau|Q^1_s-Q^2_s|^2
ds +\int_{t\wedge\tau}^\tau\|K^1_s-K^2_s\|^2ds\right)=0,
\end{eqnarray*}
for all\ $t\in[0,T]$.

It implies that for all $t\in[0,T]$
\begin{eqnarray*}
E|P^1_{t\wedge\tau}-P^2_{t\wedge\tau}|^2=0,\
E\int_0^\tau|Q^1_s-Q^2_s|^2ds=0,\ E\int_0^\tau\|K^1_s-K^2_s\|^2
ds=0.
\end{eqnarray*}

{\it Existence.}
For simplicity we assume that $f_2=0$. (In case $f_2\not=0$ we can
just smooth out $f_1$ and proceed as follows.)  Let us smooth out
$f$ to get $f^n$, i.e. let
\begin{eqnarray*}
f^n(t,P,Q,K)=\int_{\mathbb{R}^{n+n\times l}}f(t,P-n^{-1}\bar{P},Q-n^{-1}\bar{Q},K)
J(\bar{P},\bar{Q})d\bar{P}d\bar{Q},
\end{eqnarray*}
where $J(P,Q)=J_1(P)J_2(Q),$ and $J_1(P)$ is defined, for all $P\in\mathbb{R}^n$,
\begin{eqnarray*}
J_1(P)=\left\{\begin{array}{cc}
c_0\exp(-(1-|P|^2)^{-1}),\ \ &\ \ \ \mbox{as}\  |P|<1,\\
0,\ \ &\ \ \ \mbox{otherwise},
\end{array}
\right.
\end{eqnarray*}
such that the constant $c_0$ satisfies
$\int_{\mathbb{R}^{n}}J(x)dx=1$. $J_2(Q)$ is similarly defined for
any $Q\in\mathbb{R}^{n\times d}$. It easy to check that
\begin{eqnarray*}
&&|f^n(t,,P_1,Q_1,K_1)-f^n(t,P_2,Q_2,K_2)|\\
&\leq& C_n\mu\left(|P_1-P_2|+|Q_1-Q_2|+\|K_1-K_2\|\right),
\end{eqnarray*}
as $(P_i,Q_i,K_i)\in \mathbb{R}^{n}\times\mathbb{R}^{n\times
d}\times L_{\lambda(\cdot)}^2(\mathbb{R}^{n})$, $i=1,2$. Hence by
Theorem 3.3, for each $n=1,2,\cdots,$ there exists a unique solution
$(P^n_t,Q^n_t,K^n_t)$ to solve the following BDSDEP
\begin{eqnarray}\label{eq:4}
\nonumber P^n_{t\wedge\tau}&=&\xi+\int_{t\wedge\tau}^\tau
f^n(s,P^n_s,Q^n_s,K^n_s)ds+\int_{t\wedge\tau}^\tau
g(s,P^n_s,Q^n_s,K^n_s) dB_s\\
&&-\int_{t\wedge\tau}^\tau
Q^n_sdW_s-\int_{t\wedge\tau}^\tau \int_ZK^n_s(z)\widetilde{N}(dzds).
\end{eqnarray}

Applying It\^o's formula to $\left|P_t^n-P_t^m\right|^2$, we have
\begin{eqnarray*}
&
&|P_{t\wedge\tau}^n-P_{t\wedge\tau}^m|^2+\int_{t\wedge\tau}^\tau|Q_s^n-Q_s^m|^2
ds+\int_{t\wedge\tau}^\tau\|K_s^n-K_s^m\|^2ds  \\
&=&2\int_{t\wedge\tau}^\tau\left<P^n_s-P^m_s,
f^n(s,P^n_s,Q^n_s,K^n_s)-f^m(s,P^m_s,Q^m_s,K^m_s)\right>ds\\
&&
+\int_{t\wedge\tau}^\tau|g(s,P^n_s,Q^n_s,K^n_s)-g(s,P^m_s,Q^m_s,K^m_s)|^2ds\\
&&-2\int_{t\wedge\tau}^\tau\left<P^n_s-P^m_s,Q^n_s-Q^m_s\right>dW_s\\
&&+2\int_{t\wedge\tau}^\tau\left<P^n_s-P^m_s,g(s,P^n_s,Q^n_s,K^n_s)
-g(s,P^m_s,Q^m_s,K^m_s)\right>dB_s\\
&&-2\int_{t\wedge\tau}^\tau\int_Z\left<P^n_s-P^m_s,K^n_s-K^m_s\right>\widetilde{N}(
dzds)=\sum\limits_{i=1}^5I_i.
\end{eqnarray*}
Note that
\begin{eqnarray*}
I_1&=&2\int_{t\wedge\tau}^\tau\langle P^n_s-P^m_s,
\int_{\mathbb{R}^{n+n\times l}}
(f(s,P^n_s-n^{-1}\bar{P},Q^n_s-n^{-1}\bar{Q},K^n_s)\\
&&-f(s,Y^m_s-m^{-1}\bar{P},Z^m_s-m^{-1}\bar{Q},K^m_s))
J(\bar{P},\bar{Q})d\bar{P}d\bar{Q}
\rangle ds\\
&\leq& \mu\int_{t\wedge\tau}^\tau\int_{\mathbb{R}^{n+n\times l}}
((\rho(|P^n_s-P^m_s-(n^{-1}-m^{-1})\bar{P}|^2)\\
&&+|P^n_s-P^m_s-(n^{-1}-m^{-1})\bar{P}|\times(|Q^n_s-Q^m_s-(n^{-1}-m^{-1})\bar{Q}|\\
& &+\|K^n_s-K^m_s\|))+|n^{-1}-m^{-1}||\bar{P}|2)J(\bar{P},\bar{Q})
d\bar{P}d\bar{Q}ds.
\end{eqnarray*}
Since by Lemma \ref{lem:3.1} for all $n$
\begin{eqnarray*}
E\left(\sup\limits_{t\leq \tau}|P^n_t|^2+\int_0^\tau|Q^n_t|^2dt
+\int_{0}^{\tau}\|K^n_{t}\|^{2}dt\right)\leq C_T<\infty.
\end{eqnarray*}
Hence
\begin{eqnarray*}
&&E\left(|P_{t\wedge\tau}^n-P_{t\wedge\tau}^m|^2+\int_{t\wedge\tau}^\tau|Q_s^n-Q_s^m|^2ds
+\int_{t\wedge\tau}^\tau\|K_s^n-K_s^m\|^2ds\right)  \\
&\leq& \bar{C}_T(\mu^2+\mu)\int_t^T\int_{\mathbb{R}^{n+n\times l}}
(\rho(E|P^n_{s\wedge\tau}-P^m_{s\wedge\tau}-(n^{-1}-m^{-1})\bar{P}|^2)\\
& &+E|P^n_{s\wedge\tau}-P^m_{s\wedge\tau}|^2)
J(\bar{P},\bar{Q})d\bar{P}d\bar{Q} ds+\bar{C}_T(n^{-1}+m^{-1}).
\end{eqnarray*}
Note that
\begin{eqnarray*}
\rho(2E|Y^n_{s\wedge\tau}-Y^m_{s\wedge\tau}|^2+2(n^{-1}-m^{-1})^2|\bar{P}|^2)
\leq \rho(4C_T+2|\bar P|^2).
\end{eqnarray*}
But by assumption it yields that
\begin{eqnarray*}
\int\rho(4C_T+2|\bar Y|^2)J(\bar{P},\bar{Q})d\bar{P}d\bar{Q}\leq
\rho(4C_T+2)<\infty.
\end{eqnarray*}
Hence by Lemma \ref{lem:3.1} and by the Fatou lemma it is easily seen that
\begin{eqnarray*}
&&\limsup\limits_{n,m\to\infty}E|P_{t\wedge\tau}^n-P_{t\wedge\tau}^m|^2
+\limsup\limits_{n,m\to\infty}E\int_{t\wedge\tau}^\tau
\left(|Q_s^n-Q_s^m|^2+\|K_s^n-K_s^m\|^2\right)ds \\
&\leq& \widehat{C}_T(\mu^2+\mu)\int_t^T
\rho_1\left(\limsup\limits_{n,m\to\infty}
2E|P^n_{s\wedge\tau}-P^m_{s\wedge\tau}|^2\right)ds,
\end{eqnarray*}
where $\rho_1(u)=\rho(u)+u$. By the Bahari's inequality, we obtain
\begin{eqnarray*}
\limsup\limits_{n,m\to\infty}E|P_{t\wedge\tau}^n-P_{t\wedge\tau}^m|^2=0,
\ \mbox{for all}\ t\in[0,T],
\end{eqnarray*}
and
\begin{eqnarray*}
\limsup\limits_{n,m\to\infty}E\int_0^\tau
\left(|Q_s^n-Q_s^m|^2+\|K_s^n-K_s^m\|^2\right)ds=0.
\end{eqnarray*}
These, together with the Burkholder-Davis-Gundy's inequality, yield
\begin{eqnarray*}
\lim\limits_{n,m\to\infty}E\sup\limits_{0\leq t \leq
\tau}|P_{t}^n-P_{t}^m|^2=0.
\end{eqnarray*}
By the completeness of Banach space, we know that there exists a
unique $(P$, $Q,K)\in S^2([0,\tau];\mathbb{R}^{n})\times
M^2(0,\tau;\mathbb{R}^{n\times d})\times
F^{2}_{N}(0,\tau;\mathbb{R}^{n})$  such that as $n\to \infty$,
\begin{eqnarray*}
&&E\sup\limits_{0\leq t \leq \tau}|P_{t}^n-P_{t}|^2\to 0,\\
&&E\int_0^\tau|Q_s^n-Q_s|^2ds\to 0,\\
&&E\int_0^\tau \|K_s^n-K_s\|^2ds\to 0.
\end{eqnarray*}
Therefore we can take a subsequence $\{n_k\}$ of $\{n\}$, denote it
by $\{n\}$ again such that almost surely for $(t,\omega)\in
[0,T]\times \Omega$,
\begin{eqnarray*}
(P_t^n,Q_t^n,K_t^n)\to(P_t,Q_t,K_t)\ \mbox{in}\
\mathbb{R}^n\times\mathbb{R}^{n\times l}\times
L_{\lambda(\cdot)}(\mathbb{R}^n).
\end{eqnarray*}
Hence by the continuity of $f$ in $(P,Q,K)$, (H3), Lemma \ref{lem:3.1} and the
Lebesgue domination convergence theorem, we have that
\begin{eqnarray*}
E\int_0^\tau|f^n(s,P^n_s,Q^n_s,K^n_s)- f(s,P_s,Q_s,K_s)|ds\to 0,\
n\to \infty.
\end{eqnarray*}
It is easy to check that $(P,Q,K)$ is a solution of (\ref{eq:1}) by taking
the limit on both sides of (\ref{eq:4}).  \quad $\Box$

\section{Continuous dependence for solutions of BDSDEP }
\label{sec:4}

In this section, we discuss the continuous dependence for solutions of BDSDEP
(\ref{eq:1}).

\begin{theorem}\label{thm:4.1}
For $m=0,1,2,\cdots$
\begin{enumerate}
\item[(i)]
$f^m=f^m(t,p,q,k):[0,T]\times\mathbb{R}^n\times \mathbb{R}^{n\times
d}\times L^2_{\lambda(\cdot)}(\mathbb{R}^n)\to \mathbb{R}^n$ are
${\cal F}_t$-measurable such that P-a.s.
\begin{eqnarray*}
\left<p,f^m(t,p,q,k)\right>\leq \mu(t)(1+|p|^2+|p|(|q|+|k|)),
\end{eqnarray*}
where $\mu(t)$ has the property stated in (H3);
\item[(ii)]
for all $p_1,p_2\in\mathbb{R}^n; q_1,q_2\in\mathbb{R}^{n\times d};
k_1,k_2\in L^2_{\lambda(\cdot)}(\mathbb{R}^n)$, such that P-a.s.
\begin{eqnarray*}
&&\left<p_1-p_2,f^0(t,p_1,q_1,k_1)-f^0(t,p_2,q_2,k_2)\right>\\
&\leq& \mu(t)(\rho(|p_1-p_2|^2)+|p_1-p_2|(|q_1-q_2|+|k_1-k_2|)),
\end{eqnarray*}
where $\rho(\cdot)$ has the property stated in (H4);
\item[(iii)]
$g^m=g^m(t,p,q,k):[0,T]\times\mathbb{R}^n\times \mathbb{R}^{n\times
d}\times L^2_{\lambda(\cdot)}(\mathbb{R}^n)\to \mathbb{R}^{n\times
l}$ are ${\cal F}_t$-measurable such that P-a.s.
\begin{eqnarray*}
&&|g^m(t,p,q,k)|\leq \mu(t),\\
&&|g^m(t,p_1,q_1,k_1)-g^m(t,p_2,q_2,k_2)|^2\\
&\leq&\mu(t)(|p_1-p_2|^2+|p_1-p_2|(|q_1-q_2|+|k_1-k_2|)),
\end{eqnarray*}
where $\mu(t)$ has the property stated in (H3);
\item[(iv)]
$\lim\limits_{m\to \infty} \sup\limits_{p\in\mathbb{R}^n,
q\in\mathbb{R}^{n\times d},k\in L^2_{\lambda(\cdot)}(\mathbb{R}^n)}
\int_0^T|f^m(t,p,q,k)-f^0(t,p,q,k)|^2dt=0, $

$\lim\limits_{m\to \infty} \sup\limits_{p\in\mathbb{R}^n,
q\in\mathbb{R}^{n\times d},k\in L^2_{\lambda(\cdot)}(\mathbb{R}^n)}
\int_0^T|g^m(t,p,q,k)-g^0(t,p,q,k)|^2dt=0. $
\item[(v)]
 $\xi^m$ is ${\cal F}_{\tau}$-measurable and
\begin{eqnarray*}
E|\xi^m-\xi^0|^2\to 0,\ \mbox{as}\ m\to \infty,\ E|\xi^m|^2<\infty.
\end{eqnarray*}
\end{enumerate}
If $(P^m_t,Q^m_t,K^m_t)$ are solutions of the following BDSDEP: as
$0\leq s\leq T$
\begin{eqnarray*}
P^m_{s\wedge\tau}&=&\xi^m+\int_{s\wedge\tau}^{T\wedge\tau}
f^m(s,P^m_r,Q^m_r,K^m_r){\rm
d}r+\int_{s\wedge\tau}^{T\wedge\tau} g^m(s,P^m_r,Q^m_r,K^m_r)dB_r\\
&&-\int_{s\wedge\tau}^{T\wedge\tau} Q^m_r{\rm
d}W_r-\int_{s\wedge\tau}^{T\wedge\tau}
\int_ZK^m_r(z)\widetilde{N}(dzdr),\quad m=0,1,2,\cdots,
\end{eqnarray*}
then for all $0\leq s\leq T$
\begin{eqnarray*}
\lim\limits_{m\to \infty} E(\sup\limits_{s\leq r \leq T}|P^m_{r\wedge \tau} -P^0_{r\wedge
\tau}|^2+\int_{s\wedge \tau}^{T\wedge
\tau}(|Q^m_r-Q^0_r|^2+\|K^m_r-K^0_r\|^{2})dr)=0.
\end{eqnarray*}
\end{theorem}
\noindent{\bf Proof.}\ Applying It\^o's formula to $|P^m_{s} -P^0_{s}|^2$, it
follows that
\begin{eqnarray*}
&&E\left(|P^m_{s\wedge\tau}-P^0_{s\wedge\tau}|^2
+\int_{s\wedge\tau}^{T\wedge\tau}(|Q^m_r-Q^0_r|^2+\|K^m_r-K^0_r\|^2)dr\right) \\
&\leq&C_0(\mu^2(s)+\mu(s)+1)\left(\int_s^T
\rho_1(E|P^m_{s\wedge\tau}-P^0_{s\wedge\tau}|^2)dr+E|\xi^m-\xi^0|^2\right)\\
&&+C_0E\int_s^T|f^m(P^m_r,Q^m_r,K^m_r)-f^0(P^m_r,Q^m_r,K^m_r)|^2dr\\
&&+C_0E\int_s^T|g^m(P^m_r,Q^m_r,K^m_r)-g^0(P^m_r,Q^m_r,K^m_r)|^2dr,
\end{eqnarray*}
where
\begin{eqnarray*}
\rho_1(u)=\rho(u)+u.
\end{eqnarray*}
Hence
\begin{eqnarray*}
\limsup\limits_{m\to \infty}E|P^m_{s\wedge\tau}-P^0_{s\wedge\tau}|^2
\leq C_0\int_s^T\hat \mu(r)\rho_1\left(\limsup\limits_{m\to
\infty}E|P^m_{r\wedge\tau}-P^0_{r\wedge\tau}|^2\right)dr.
\end{eqnarray*}
By the Bahari's inequality, we obtain
\begin{eqnarray*}
\limsup\limits_{m\to
\infty}E|P^m_{s\wedge\tau}-P^0_{s\wedge\tau}|^2=0.
\end{eqnarray*}
It is easily derived that
\begin{eqnarray*}
\lim\limits_{m\to \infty} E\left(\int_{s\wedge \tau}^{T\wedge
\tau}(|Q^m_r-Q^0_r|^2+\|K^m_r-K^0_r\|^{2})dr\right)=0.
\end{eqnarray*}
Applying It\^o's formula to $|P^m_{s} -P^0_{s}|^2$ on $[0,t\wedge \tau]$, we have
\begin{eqnarray*}
&&E\sup\limits_{t \leq \tau}|P^m_{t\wedge \tau} -P^0_{t\wedge \tau}|^2\\
&\leq& E|P^m_{0} -P^0_{0}|^2+\int_0^{\tau}|g^m(s,P^m_s,Q^m_s,K^m_s)-g^0(s,P^0_s,Q^0_s,K^0_s)|^2ds\\
&&+2E\int_0^{\tau}\langle
P^m_{s} -P^0_{s},f^m(s,P^m_s,Q^m_s,K^m_s)-f^0(s,P^0_s,Q^0_s,K^0_s)\rangle ds\\
&&-\int_0^{\tau}|Q^m_{s} -Q^0_{s}|^2ds +\int_0^{\tau}\|K^m_{s}-K^0_{s}\|^2ds\\
&&+2E\sup\limits_{t \leq \tau}\int_0^{t\wedge
\tau}\langle P^m_{s} -P^0_{s},g^m(s,P^m_s,Q^m_s,K^m_s)-g^0(s,P^0_s,Q^0_s,K^0_s)\rangle dB_s\\
&&-2E\sup\limits_{t \leq \tau}\int_0^{t\wedge\tau}\langle P^m_{s} -P^0_{s},Q^m_{s} -Q^0_{s}\rangle
dW_s\\
&&+2E\sup\limits_{t \leq \tau}\int_0^{t\wedge\tau}\int_Z\langle
P^m_{s} -P^0_{s},K^m_s(z)-K^0_s(z)\rangle\widetilde{N}(dzds).
\end{eqnarray*}
By the similar arguments in Lemma \ref{lem:3.1}, we obtain
\begin{eqnarray*}
\lim\limits_{m\to \infty} E\sup\limits_{s\leq r \leq T}|P^m_{r\wedge \tau} -P^0_{r\wedge
\tau}|^2=0.
\end{eqnarray*}
\quad $\Box$

We also have the other useful continuous dependence
for solutions of BDSDEP as follows:

\begin{theorem}\label{thm:4.2}
For $m=0,1,2,\cdots$
\begin{enumerate}
\item[(i)]
$f^m(t,p,q,k)$ are ${\cal F}_t$-measurable such that P-a.s.
\begin{eqnarray*}
|f^m(t,p,q,k)|\leq C_0(1+|p|+|q|+|k|),
\end{eqnarray*}
where $C_0\leq 0$ is a constant;
\item[(ii)]
for all $p_1,p_2\in\mathbb{R}^n; q_1,q_2\in\mathbb{R}^{n\times l};
k_1,k_2\in\mathbb{R}^n$, such that
\begin{eqnarray*}
&&\left<p_1-p_2,f^m(t,p_1,q_1,k_1)-f^m(t,p_2,q_2,k_2)\right>\\
& &\leq
\mu(t)\rho(|p_1-p_2|^2)+|p_1-p_2|^2(|q_1-q_2|^2+|k_1-k_2|^2),
\end{eqnarray*}
where $\mu(t)$ has the property stated in (H3)
and $\rho(\cdot)$ has tthe property stated in (H4);
\item[(iii)]
 The same as (iii) in Theorem \ref{thm:4.1};
\item[(iv)]
$\lim\limits_{m\to \infty}f^m(t,p,q,k)=f^0(t,p,q,k)$, P-a.s.

$\lim\limits_{m\to \infty}g^m(t,p,q,k)=g^0(t,p,q,k)$, P-a.s.
\item[(v)]
 The same as (v) in Theorem \ref{thm:4.1};
\end{enumerate}
Then the conclusion of Theorem \ref{thm:4.1} still holds.
\end{theorem}

The proof can be completed similarly as that of Theorem \ref{thm:4.1}.

\section{The probabilistic interpretation of SPDIEs}
\label{sec:5}

The connection of BDSDEs and systems of second-order quasilinear
SPDEs was observed by Pardoux and Peng \cite{PP2}. This can be regarded as
a stochastic version of the well-known Feynman-Kac formula which
gives a probabilistic interpretation for second-order SPDEs of
parabolic types. Thereafter this subject has attracted many
mathematicians, referred to Bally and Matoussi \cite{BM}, Zhang and Zhao
\cite{ZZ}, Hu and Ren \cite{HR}, see also Ren et al. \cite{RLH}. In \cite{HR}, the authors
got a probabilistic interpretation for the solution of a semilinear
SPDIE, via BDSDEs with L$\acute{\rm e}$vy process for a fixed
terminal time under Lipschitzian assumption. This section can
be viewed as a continuation of such a theme, and will exploit the
above theory of BDSDEP with non-Lipschitzian coefficients and random
terminal time in order to provide a probabilistic formula for the
solution of a quasilinear SPDIE.

Let $D$ be a bound domain in $\mathbb{R}^m$ with boundary $\partial
D=S$.

First, consider the following forward SDE with Poisson jumps in
$\mathbb{R}^m$ for any given $(t,x)\in [0,T]\times D$
\begin{eqnarray}\label{eq:5}
\nonumber X_t &=& x+\int_t^s b(r,X_r)dr+\int_t^s g(r,X_r)dW_r\\
&&+\int_t^s \int_Z
h(r_-,X_{r_-},z)\widetilde{N}(dz{\rm d}r),\quad t\leq s\leq T,
\end{eqnarray}
where
\begin{eqnarray*}
b:[0,T]\times \mathbb{R}^m \to \mathbb{R}^m,\ \sigma:[0,T]\times
\mathbb{R}^m \to  \mathbb{R}^{m\times d},\ h:[0,T]\times
\mathbb{R}^m\times Z \to \mathbb{R}^{m}.
\end{eqnarray*}
It is known that, if coefficients are less than linear increasing,
and satisfy the Lipschitz condition, then SDE (\ref{eq:5}) has a unique
solution. (See \cite{S})

Now for any $(t,x)\in [0,T]\times D$, let
\begin{eqnarray*}
\tau=\tau_x=\inf\{s>t:X_s^{t,x} \notin D\},\ \mbox{and}\
\tau=\tau_x=T,\ \mbox{for}\ \inf\{\phi\}.
\end{eqnarray*}
Consider the following BDSDEP (for simplicity, denote
$X_s=X_s^{t,x}$),
\begin{eqnarray}\label{eq:6}
\nonumber P_s&=&\Phi(X_\tau)+\int_{s\wedge
\tau}^{\tau}f(r,X_r,P_r,Q_r,K_r) dr+\int_{s\wedge
\tau}^{\tau}g(r,X_r,P_r,Q_r,K_r)dB_r\\
&&+\int_{s\wedge \tau}^{\tau}Q_rdW_r+\int_{s\wedge
\tau}^{\tau}\int_Z K_{r_-}(z)\widetilde{N}(dzdr),\quad t\leq s\leq
T,
\end{eqnarray}
where
\begin{eqnarray*}
&&f:[0,T]\times \mathbb{R}^m\times \mathbb{R}^n\times
 \mathbb{R}^{n\times d}\times L^2_{\lambda(\cdot)}(
\mathbb{R}^n) \to \mathbb{R}^n,\\
&&g :[0,T]\times \mathbb{R}^m\times\mathbb{R}^n\times
 \mathbb{R}^{n\times d}\times L^2_{\lambda(\cdot)}(
\mathbb{R}^n)\to \mathbb{R}^{n\times l},\\ && \Phi:\mathbb{R}^m \to
\mathbb{R}^n.
\end{eqnarray*}
Suppose that $f(t,x,\cdot,\cdot,\cdot)$ and
$g(t,x,\cdot,\cdot,\cdot)$ satisfy the conditions in Theorem \ref{thm:3.4}
uniformly for $t$ and $x$, and suppose that
$E|\Phi(X_\tau)|^2<\infty$, then by Theorem \ref{thm:3.4}, BDSDEP (\ref{eq:6}) has a
unique solution $(P_t,Q_t,K_t)\in S^2([0,\tau];\mathbb{R}^{n})\times
M^2(0,\tau;\mathbb{R}^{n\times d})$ $\times
F^{2}_{N}(0,\tau;\mathbb{R}^{n})$.

We now relate BDSDEP (\ref{eq:6}) to the following system of quasilinear
second-order parabolic SPDIE:
\begin{eqnarray}\label{eq:7}
\left\{
\begin{array}{lllll}
{\cal L}u(t,x)dt\\
=f(t,x,u(t,x),\nabla u(t,x)g(t,x),u(t,x+h(t,x,\cdot))-u(t,x))dt\\
+g(t,x,u(t,x),\nabla
u(t,x)g(t,x),u(t,x+h(t,x,\cdot))-u(t,x))dB_t,\\
\forall (t,x)\in [0,T]\times D\\
u(T,x)=\Phi(x),\quad \forall (t,x)\in [0,T]\times D;\\
u(t,x)=\Phi(x),\quad \forall (t,x)\in [0,T]\times S,\\
\end{array}
\right.
\end{eqnarray}
where $u:\mathbb{R}_+\times\mathbb{R}^m\to \mathbb{R}^n$,
\begin{eqnarray*}
{\cal L}u =\left(\begin{array}{c} Lu_1\\\vdots\\
Lu_n\end{array}\right),
\end{eqnarray*}
with
\begin{eqnarray*}
&&Lu_k(t,x)\\
&=&\displaystyle\frac{\partial u_k}{\partial
t}(t,x)+\sum\limits_{i=1}^nb_{i}(t,x) \displaystyle\frac{\partial
u_k }{\partial x_i}(t,x)+
\displaystyle\frac{1}{2}\sum\limits_{i,j=1}^n
(\sigma\sigma*)_{ij}(t,x)\displaystyle\frac{\partial^2 u_k
}{\partial x_i\partial
x_j}(t,x)\\
&&+\displaystyle\int_Z(u_k(t,x+h(t,x,z))-u_k(t,x)-\sum\limits_{i=1}^n
h_i(t,x,z)\displaystyle\frac{\partial u_k}{\partial x_i}(t,x)
)\lambda(dz),\\
&&k=1,\cdots,n.
\end{eqnarray*}
Now assume that $\sigma$ is uniformly non-degenerate, i.e. there
exists a constant $\beta > 0$, such that
\begin{eqnarray*}
\displaystyle\frac{1}{2}\sum\limits_{i,j=1}^n
(\sigma\sigma*)_{ij}(t,x)\xi_i\xi_j\geq \beta|\xi|^2,\ \mbox{for
all}\ \xi\in\mathbb{R}^m,\ \mbox{and}\  (t,x)\in [0,T]\times \bar D,
\end{eqnarray*}
where $\bar D=$ the closure of $D$. Hence, SPDIE (\ref{eq:7}) is a true
quasilinear type equation.

We can assert that

\begin{theorem}\label{thm:5.1}
Under the above related conditions, and $b$, $\sigma$, $h$, $f$ and
$g$ are of class $C^3$, and $\Phi$ is of class $C^2$. Suppose SPDIE
(\ref{eq:7}) has a unique solution $u(t,x)\in C^{1,2}(\Omega\times
[0,T]\times \mathbb{R}^m; \mathbb{R}^n)$. Then, for any given
$(t,x)$, $u(t,x)$ has the following interpretation
\begin{eqnarray}\label{eq:8}
u(t,x)= P_t,\end{eqnarray} where $P_t$ is determined uniquely by
(\ref{eq:5}) and $(\ref{eq:6})$.
\end{theorem}

\noindent{\bf Proof}\ Applying It\^o's formula to $u(t,X_t)$
(see Theorem 6 in \cite{S}) on $[s\wedge \tau, \tau]$, we obtain
\begin{eqnarray*}
&&u(\tau,X_\tau)-u(s\wedge \tau,X_{s\wedge \tau})\\
&=&\int_{s\wedge\tau}^\tau\frac{\partial u}{\partial r}(r,X_r)dr+\int_{s\wedge
\tau}^\tau\sum\limits_{i=1}^mb_{i}(r,X_r)
\displaystyle\frac{\partial u}{\partial x_i}(r,X_r)dr
\\
&&+\int_{s\wedge \tau}^\tau\nabla
u(r,X_r)\sigma(r,X_r)dW_r+\int_{s\wedge
\tau}^\tau\displaystyle\frac{1}{2}\sum\limits_{i,j=1}^m
(\sigma\sigma*)_{ij}(r,X_r)\displaystyle\frac{\partial^2u }{\partial
x_i\partial x_j}(r,X_r)dr\\
&& +\int_{s\wedge
\tau}^\tau\int_Z(u(r,X_r+h(r,X_r,z))-u(r,X_r))\widetilde{N}(dzdr)\\
&&+\int_{s\wedge
\tau}^\tau\displaystyle\int_Z(u(r,X_r+h(r,X_r,z))-u(r,X_r)\\
&&-\sum\limits_{i=1}^m
h_i(s,X_s,z)\displaystyle\frac{\partial u}{\partial x_i}
(r,X_r))\lambda(dz)dr.
\end{eqnarray*}
Because $u(t,x)$ satisfies SPDIE (\ref{eq:7}), it holds that
\begin{eqnarray*}
&&\Phi(X_\tau)-u(s\wedge \tau,X_{s\wedge \tau})\\
&=&\int_{s\wedge\tau}^\tau f(r,X_r,u(r,X_r),\nabla
u(r,X_r)\sigma(r,X_r,u(r,X_r)),\\
&&u(r,X_r+h(r,X_r,\cdot))-u(r,X_r))dr\\
&&+ \int_{s\wedge \tau}^\tau g(r,X_r,u(r,X_r),\nabla
u(r,X_r)\sigma(r,X_r,u(r,X_r)),\\
&&u(r,X_r+h(r,X_r,\cdot))-u(r,X_r))dB_r
+\int_{s\wedge \tau}^\tau\nabla u(r,X_r)g(r,X_r)dW_r\\
&& +\int_{s\wedge
\tau}^\tau\displaystyle\int_Z(u(r,X_r+h(r,X_r,z))-u(r,X_r))\widetilde{N}(dzdr).
\end{eqnarray*}
It is easy to check that $(u(t,X_t),\nabla u(t,X_t)\sigma(t,X_t),
u(t,X_t+h(t,X_t,\cdot))-u(t,X_t))$ coincides with the unique
solution of BDSDEP (\ref{eq:6}). It follows that
\begin{eqnarray*}
u(t,x)= P_t.
\end{eqnarray*}
\quad $\Box$

\begin{remark}\label{rmk:5.2}
(\ref{eq:8})can be called a stochastic Feynman-Kac formula for SPDIE (\ref{eq:7}), which is
a useful tool in the study of SPDIE.\end{remark}

\end{document}